\documentclass{article}

\usepackage{amsfonts}

\usepackage{amsmath}

\usepackage[backend=biber,style=numeric]{biblatex}
\usepackage[english]{babel}

\begin{document}

\centerline{\Large\bf The Tate Thomason Conjecture}

\vspace{20pt}

\centerline{\bf Marcelo Gomez Morteo}

\vspace{20pt}

\footnote{ email valmont8ar@hotmail.com}

\pagestyle{myheadings}

\vspace{16pt}

\begin{abstract}

\vspace{16pt}

We prove the Tate Thomason conjecture using $\cal K$ where $\cal K$ is the complex topology spectrum. ( see the introduction below). Fundamental to our proof is Theorem 2.2. Much of this work is related to the article [7]

\vspace{16pt}


\vspace{6pt}

\end{abstract}

\vspace{16pt}

\textbf{Introduction:}

\vspace{16pt}

By a spectrum we mean the following: A spectrum X is a collection of simplicial sets $X_n$ for $n\geq0$ together with morphisms of simplicial sets $\sigma\sb{n} :\Sigma X\sb{n} \mapsto Y\sb{n}$. A morphism of spectra $f:X \mapsto Y$ is a collection of morphisms $f_n:X_n \mapsto Y_n$ of simplicial sets that commute with the structure maps $\sigma_n$, ie $\sigma_n \circ \Sigma f_n=f\sb{n+1} \circ \sigma_n$. We are going to consider the stable homotopy category $\cal S$ derived from the Bousfield-Friedlander model category of simplicial spectra.
The topological K spectrum is written $\cal K$ and its localization is written $L$ most of the time.
The spectra we consider are K(X) where X is a scheme over a finite field $F_q$ with $q=p^s$ where $p\neq l$, and K(X) is Quillen's Algebraic K theory spectra related to X. See [3] pages 148-150.
Our aim is to prove the Tate Thomason conjecture ([7]): Let $X\sb{\infty}$ be $X\times\sb{Spec(F_q)}Spec(\bar{F_{q}})$, where X is a smooth projective variety over $Spec(F_q)$, and $\bar{F_{q}}$ is the algebraic closure of $F_{q}$, then the homotopy group $\pi\sb{-1}(LK(X\sb{\infty}))$ is reduced. This statement implies the Tate Conjecture, ( See Remark $0.2$ below and also see [7] page 390 diagram 21). The proof is based on theorem 2.2  from section 2.

\vspace{30pt}
\textbf{0. The Tate Conjecture}

\vspace{16pt}

For the projective smooth variety $X\sb{\infty}$ over the finite field $F_{q}$, consider the $i$-th Chow Group $CH^{i}(X\sb{\infty})$  generated by cycles of codimension $i$ on $X\sb{\infty}$ modulo rational equivalence. If $l$ is a prime invertible in $F_{q}$ we can define a $Q_{l}$-linear cycle morphism from the $i$-th Chow group into the $2i$-th $l$-adic etale cohomology with $i$-th Tate twist coefficient:

\vspace{16pt}

\[\gamma^{i}_{Q_{l}}:CH^{i}(X\sb{\infty})\otimes Q_{l} \mapsto H^{2i}(X\sb{\infty},Q_{l}(i))\]

\vspace{16pt}

 Consider $F_{q}$ as a finite subfield of $\bar{F_{q}}$. Let $H^{2i}(X\sb{\infty},Q_{l}(i))^{Gal(\bar{F_{q}}/F{q})}$ be the subspace of $H^{2i}(X\sb{\infty}, Q_{l}(i))$ of finite order  under the action of $Gal(\bar{F_{q}}/F_{q})$, then the Tate Conjecture states:

\vspace{16pt}

0.1 \textbf{Tate Conjecture}: \emph{The image of the cycle class morphism $\gamma^{i}_{Q_{l}}$ in the cohomology group $H^{2i}(X\sb{\infty}, Q_{l}(i))$ is exactly
 $H^{2i}(X\sb{\infty}, Q_{l}(i))^{Gal(\bar{F_{q}}/F_{q})}$}

\vspace{30pt}

0.2 \emph{Thomason's reformulation of the Tate Conjecture:}

\vspace{16pt}

Take the smooth projective variety $X_{\infty}$. Let $ X_{n}=X\otimes F_{q^{n}}$ where $X$ is a smooth projective variety over $F_{q}$ as mentioned in the introduction.   Let $\cal K$ be the complex topology spectrum. Thomason proves in [7] the following theorem:

\vspace{16pt}

\textbf{Theorem 0.1 (Tate Thomason's Conjecture):} \emph{The Tate Conjecture  is equivalent to the finiteness statement that for all $n\in N$,  See [7]}

\vspace{20pt}

\[Hom(Q/Z_{(l)}, \pi_{-1}(L_{\cal K}K(X_{n})))=0\]

\vspace{10pt}

Remark 0.1: observe that in [7] page 390 diagram 21 Thomason states that the arrow $f\otimes Q$ should be $0$ for the Tate conjecture to be true. Hence we only need to prove that

\vspace{10pt}

\[Hom(Q/Z_{(l)}, \pi_{-1}(L_{\cal K}K(X\sb{\infty})))=0\]
\vspace{10pt}

Remark 0.2: Though the above Hom being $0$ is sufficient to prove the Tate Conjecture, to prove the equivalence with the Tate Conjecture we need the Hom with the $X_{n}$ See [7]

\vspace{10pt}

Definition 0.1: \emph{A group $G$ which verifies that $ Hom (Q/Z_{(l)}, G)=0$ is said to be $l$ reduced.}

\vspace{10pt}

We will simplify all over this work the terminology by saying that a group is reduced when we really mean that it is $l$-reduced.

\vspace{20pt}

The main lemmas stated by Thomason in [7] to prove  the Tate Conjecture are

\vspace{20pt}

Lemma 0.1: \emph{Let $X$ be a smooth projective variety over $F_{q}$ with $q=p^{n}$ and $p$ prime, or $X_{\infty}$ If $l$ is a prime number different from $p$, and if $K^{Top}(X)$ is in Thomason's notation the topological $K$-theory spectrum, and $(...)\widehat{}$ is the $l$-adic completion of a spectrum,
then}

\vspace{16pt}

\[L_{\cal K}K(X)^{\widehat{}}\cong K^{Top}(X)^{\widehat{}}\]

\vspace{30pt}

The proof of this lemma follows from Thomason's descent theorem proved in [6, 4.1]which relates algebraic K- theory to topological K-theory. See also ([7], page 388, equation (14) and references therein)

\vspace{20pt}

Remark 0.3: It also follows from [7] page 388 equation (13) that

\vspace{16pt}

\[K^{Top}/l^{\nu}(X) \cong L_{\cal K}(X)/l^{\nu}\]

\vspace{10pt}

See also the descent problem in [4] section 4)

\vspace{30pt}

Lemma 0.2:( See [7],the lemma of page 387)

\vspace{5pt}

\emph{The image of $colim K_{0}^{Top}(X_{n})^{\widehat{}} \otimes Q$ in $K_{0}^{Top}(X_{\infty})^{\widehat{}} \otimes Q$ consists precisely of those elements of finite orbit under $Gal (\bar{F_{q}}/F_{q})$.}

\vspace{30pt}

Remark 0.4:(See [7], page 387 square diagram 10) which is equal to

\vspace{16pt}

\[
colim (K_{0}(X_{n})\otimes Q_{l})=K_{0}(X_{\infty})\otimes Q_{l} \rightarrow colim  K_{0}^{Top}(X_{n})^{\widehat{}}\otimes Q) \downarrow K_{0}^{Top}(X_{\infty})^{\widehat{}}\otimes Q \]

\vspace{20pt}

Where the equality means an isomorphism and in light of this, lemma 0.2, and the fact that $K_{0}(X_{\infty})\otimes Q_{l}$ is isomorphic to $\bigoplus_{i=1}^{d}CH^{i}(X_{\infty})\otimes Q_{l}$ with $d$ the dimension of $X$, and that $K_{0}^{Top}(X_{\infty})^{\widehat{}}\otimes Q$ is isomorphic to $ \bigoplus_{i=1}^{d}H_{et}^{2i}(X_{\infty},Q_{l}(i))$,([7]) the above composed map is exactly the cycle map $\gamma(X_{\infty})_{Q_{l}}$ and the Tate Conjecture is equivalent to the conjecture that $K_{0}(X_{\infty})\otimes Q_{l}$ and
$colim (K_{0}^{Top}(X_{n})^{\widehat{}}\otimes Q)$ have the same image in $K_{0}^{Top}(X_{\infty})^{\widehat{}}\otimes Q$.

\vspace{16pt}
Lemma 0.3: is diagram 21 on page 390 of [7]

\vspace{20pt}

Remark 0.5: To our knowledge, Thomason never said anything about the $l$-reducibility of the homotopy group $\pi_{-1}L_{\cal K}K(X_{\infty}))$

\vspace{30pt}

\textbf{1. The Moore Spectrum and an important exact sequence }

\vspace{20pt}

\emph{Remark 1.1:} Given $l^{\nu}$, with $\nu\in N$,and $M(Z/l^{\nu})$ the Moore spectrum of the ring $Z/l^{\nu}$,and if $X$ is a smooth variety over a field $k$ where $k=F_{q^{n}}$ or $k=\bar{F_{q}}$ with $q=p^{n}$ and L is from now on the localization functor at the complex K-Theory spectrum  $\cal K$, then

\vspace{12pt}

$LK(X)/{l^{\nu}}=LK(X)\wedge M(Z/{l^{\nu}})=K(X)\wedge L\Sigma^{\infty} S^{o}\wedge M(Z_/{l^{\nu}})=K(X)/{l^{\nu}}\wedge L\Sigma^{\infty} S^{o}=L(K(X)/{l^{\nu}})$, equalities following from the smash property for L.

\vspace{12pt}

\emph{Remark 1.2:} $\pi_{*}(L(K(X)/{l^{\nu}}))$ is $l^{\nu}$-torsion.

\vspace{12pt}

The claim follows from the exact sequence

\vspace{4pt}

\[0 \mapsto \pi_{*}(LK(X)\otimes Z/l^{\nu} \mapsto \pi_{*}(LK(X)\wedge M(Z/l^{\nu})) \mapsto \]

\vspace{2pt}

\[ \mapsto Tor^{1} (\pi_{*-1}(LK(X)),Z/l^{\nu}) \mapsto 0\]

\vspace{10pt}
which splits (See [6]  Appendix A, (6) )

\vspace{16pt}

2. \textbf{Main Theorem.}

\vspace{16pt}

Theorem 2.2: \emph{ The Tate Module ( which we call $M$) of the abelian group $\pi_{m}(LK(X_{\infty}))$ is trivial for every integer $m$}

\vspace{10pt}
As we have mentioned above, we have by [6] Appendix 6) that the exact sequence below with $\pi_{m}(LK(X_{\infty}))\wedge M(Z/l^{\nu})=\pi_{m}(LK(X_{\infty})/l^{\nu})$  splits

\vspace{10pt}

\[ 0 \mapsto \pi_{m}(LK(X_{\infty}))\otimes Z/l^{\nu} \mapsto \pi_{m}(LK(X_{\infty})/l^{\nu}) \mapsto \tag{1} \]

\vspace{2pt}

\[\mapsto Tor^{1} (\pi_{m-1}(LK(X_{\infty})),Z/l^{\nu}) \mapsto 0 \]

 \vspace{10pt}

 Taking inverse limit over $\nu$ on the above exact sequence, we get an exact sequence since $lim$ is left exact, and moreover exact in the case of the above exact sequence since the left extreme of that sequence is a surjective inverse system, and therefore Mittag-Leffler. Then, if we prove that

 \vspace{10pt}

\[lim \pi_{m}(LK(X_{\infty}))\otimes Z/l^{\nu}=(\pi_{m}(LK(X_{\infty})))^{l}\cong lim \pi_{m}(LK(X_{\infty}))\wedge M(Z/l^{\nu}) \]

 \vspace{10pt}

 where $(\pi_{m}(LK(X_{\infty})))^{l}$ stands for $l$-adic completion of that abelian group, with $l$ prime to $p$ the characteristic of the base field of our smooth projective variety, being the $lim$ exact here, we conclude that the Tate Module of
 $\pi_{m}(LK(X_{\infty}))$ is trivial since it is equal, with inverse limit over $\nu$ for fixed $l$ to

 \vspace{5pt}

 \[lim Tor^{1} (\pi_{m-1}(LK(X_{\infty})),Z/l^{\nu})\]

 \vspace{10pt}

Again,taking inverse limit with respect to the parameter $\nu$ with $l$ fixed,

\vspace{8pt}

\[lim Tor^{1} (\pi_{m-1}LK(X_{\infty})),Z/l^{\nu}) )=M={0}\].

\vspace{8pt}

where $ M=\lim Tor^{1} (\pi_{m-1}LK(X_{\infty}),Z/l^{\nu})$ is equal to
\[{\Pi_{i=1}^{\infty}\{{ g_{i}=l^{i}-torsion-element\in \pi_{m-1}(LK(X_{\infty}))}: lg_{i+1}=g_{i}}\}\]

\vspace{8pt}
by definition of the inverse limit and since $Tor^{1}(-, Z/l^{\nu})$ are the $l^{\nu}$ -torsion
elements of $ \pi_{m-1}LK(X_{\infty})$

\vspace{10pt}

(See [5]  8.4 pp 223 and also that $Tor_{Z}^{1}(A,B)=Tor_{Z}^{1}(B,A)$ on pp 222).

\vspace{12pt}

 Now, consider the image of a map $f$ in $Hom(Q/Z_{(l)},\pi_{m-1}(LK(X_{\infty}))$

 \vspace{12pt}

 Since $Q/Z_{(l)}$ is an $l$ quasicycle divisible group , if $D$ is the image of $f$ then $D$ is an increasing union of groups $D_{n}$, of order $l^{n}$, where each $D_n$ is generated by one generator $g_{n}$ wich is $l^{n}$-torsion, and also $lg_{n}=g_{n-1}$, if the Tate module $M$ trivial, this sequence of $g_{n}$ has to be trivial and therefore $f=0$

 \vspace{8pt}

Therefore details of the Tate Module of the homotopy groups $\pi_{m-1}LK(X_{\infty})$ are of fundamental importance. The above arguments show that the Tate module $M$ is isomorphic to $Hom(Q/Z_{(l)},\pi_{m-1}(LK(X_{\infty}))$ by taking the sequence of generators of the increasing sequence of cyclic subgroups of $Q/Z_{(l)}$ to the sequence of it's images which belong to $M$. By [7] pp 389 (15) we have an exact sequence

 \vspace{10pt}

\[0 \mapsto Ext(Q/Z_{(l)},\pi_{m}(LK(X_{\infty}))) \mapsto \pi_{m}((LK(X_{\infty}))^{l})\mapsto \tag{2} \]

 \vspace{2pt}

\[ \mapsto Hom(Q/Z_{(l)},\pi_{m-1}(LK(X_{\infty}))\mapsto 0 \]

\vspace{10pt}

We can compare the exact sequence (2) with the above exact sequence (1) after taking $lim$ ( which is equation (3)) and using the fact that $lim \pi_{m}(LK(X_{\infty}))\wedge M(Z/l^{\nu})\cong \pi_{m}((LK(X_{\infty}))^{l})$ by [7] (6) pp 387 and (13) pp 388 obtaining

\vspace{5pt}

\[0\mapsto (\pi_{m}(LK(X_{\infty})))^{l} \mapsto \pi_{m}((LK(X_{\infty}))^{l})\mapsto \tag{3} \]

\vspace{3pt}

\[\mapsto lim Tor^{1} (\pi_{m-1}LK(X_{\infty}),Z/l^{\nu})\mapsto 0 \]

\vspace{10pt}

We have the same middle term in both exact sequences (2) and (3) and moreover we also proved that the extreme right terms are isomorphic. Therefore by the snake lema, the extreme left terms must also be isomorphic. By [7] pp 389 (16) we have another exact sequence:

\vspace{8pt}

\[0\mapsto lim^{1}Hom(Z/l^{k}Z,\pi_{m}((LK(X_{\infty}))))\mapsto Ext(Q/Z_{(l)},\pi_{m}(LK(X_{\infty})))\mapsto \tag{4} \]

\vspace{2pt}

\[lim \pi_{m}(LK(X_{\infty}))\otimes Z/l^{\nu} \mapsto 0 \]

\vspace{8pt}
Since we have just shown by comparing the exact sequences (2) and (3) that $Ext(Q/Z_{(l)},\pi_{m}(LK(X_{\infty})))$ , is isomorphic to the right extreme of the above exact sequence (4), then in (4) the left extreme must be trivial, ie

\vspace{8pt}

\[lim^{1}Hom(Z/l^{k}Z,\pi_{m}(LK(X_{\infty})))=0 \tag{5} \]

\vspace{8pt}
 Thomason states in [7] that if the $l$ torsion of the groups $\pi_{m}(LK(X_{\infty}))$ are finite, then the Mittag Leffler condition holds for the inverse system of the abelian groups $Hom(Z/l^{k}Z,\pi_{m}((LK(X_{\infty}))))$ and then $lim^{1}=0$ for those groups. It is a well known fact that the Mittag Leffler condition for an inverse system of abelian groups implies the statement $lim^{1}=0$ for that system of abelian groups and that is what is implied in Thomason statement in [7]. The other way around is not always true, since an inverse system having $lim^{1}=0$ does not in general imply the Mittag Leffler condition for that inverse system. But if the abelian groups  $\pi_{m}((LK(X_{\infty})))$ are at most countable abelian groups then $Hom(Z/l^{k}Z,\pi_{m}((LK(X_{\infty}))))$ is an inverse system of at most countable abelian groups. In [1] it is stated,  and in [2] it is proven, that an inverse system of at most countable abelian groups has $Lim^{1}=0$  if and only if the Mittag Leffler condition holds true for that inverse system. Using this last statement if the $Hom(Z/l^{k}Z,\pi_{m}((LK(X_{\infty}))))$ are at most countable, since  $lim^{1}Hom(Z/l^{k}Z,\pi_{m}((LK(X_{\infty}))))=0$ then the Mittag Leffler must hold true for the inverse system $Hom(Z/l^{k}Z,\pi_{m}((LK(X_{\infty}))))$  in which case the homotopy groups $\pi_{m}(LK(X_{\infty}))$  must have bounded $l$ torsion, since the inverse system of our hom groups, being Mittag Leffler at certain point become stationary.If that is the case, ie if they have bounded $l$-torsion, taking an $f$ in $Hom(Q/Z_{(l)},\pi_{m-1}(LK(X_{\infty})))$ it's image is $l$-torsion, and hence must have bounded $l$-torsion because of what we have just proved. But that is not possible since $Q/Z_{(l)}$ is divisible and therefore the image of $f$ is an $l$- torsion and divisible abelian group  .Hence $f=0$ and we conclude that the Tate module M is trivial as stated in theorem 2.2, and moreover we can state the following theorem:

\vspace{30pt}

 Theorem 2.3: \emph{$Hom(Q/Z_{(l)},\pi_{m-1}(LK(X_{\infty}))=0$  for all $m$}

 \vspace{20pt}

   To end the proof of theorem 2.2 and theorem 2.3 ,we must show that the inverse system of hom groups $Hom(Z/l^{k}Z,\pi_{m}((LK(X_{\infty}))))$ is an inverse system of at most countable abelian groups.The image of an $f$ belonging to the abelian group $Hom(Q/Z_{(l)},\pi_{m}(LK(X_{\infty}))$ can only intersect $l$-quasicyclic subgroups of $\pi_{m}((LK(X_{\infty})))$ so that concerning the $l$ quasicyclyc part of $\pi_{m}((LK(X_{\infty})))$,this fact informs us that $\pi_{m}((LK(X_{\infty})))$ must have at least an $l$ quasicyclic direct summand as a subgroup. If $Hom(Z/l^{k}Z,\pi_{m}((LK(X_{\infty}))))$ is uncountable then writing this homotopy groups as $A_{r}\times A_{d}$ where $A_{r}$ is the reduced subgroup and $A_{d}$ the divisible subgroup of $\pi_{m}((LK(X_{\infty})))$  then either b) $Hom(Z/l^{k}Z,A_{r})$ is uncountable or a) $Hom(Z/l^{k}Z,A_{d})$ is uncountable. In case a) the $l$ quasicyclic part of $\pi_{m}((LK(X_{\infty})))$ which is included in $A_{d}$  must have an uncountable number of direct summands or must be an  infinite product of infinite $l$ quasicyclic groups. We will focus first on case a), and later on we will consider case b). In  case a) we must analyze if the $l$ quasicyclic part of $\pi_{m}((LK(X_{\infty})))$ can be i) an infinite product of $l$ quisicyclic groups or ii) an uncountable direct sum of $l$ quasicyclic groups. We might as well consider as $\pi_{m}((LK(X_{\infty})))$ having one of those two structures, since adding a complementary abelian group will basically not change what we will prove in what follows because since the exact sequence given below splits by [6] appendix (6) (something we have already mentioned above) then the map from the middle term to the right hand side term of that exact sequence is a retraction, or equivalently the right hand side term is mapped by an injective morphism to its middle term something which is not possible since as we will soon point out the middle term is a finite abelian group while the right hand side term cannot be finite either by adding or not adding to an infinite product of $l$ quasicyclic groups a complementary abelian group. Also, since both cases are similar we analyze only  case i). Hence,if we consider the case of an infinite product of $l$-quasicyclic groups, we are going to show that that structure cannot possibly hold. As just stated above, reconsider once more the already mentioned exact sequence (1)  below theorem 2

 \vspace{30pt}

\[0 \mapsto \pi_{m}(LK(X_{\infty}))\otimes Z/l^{\nu} \mapsto \pi_{m}(LK(X_{\infty})/l^{\nu}) \mapsto \tag{1} \]

\vspace{3pt}

  \[Tor^{1} (\pi_{m-1}(LK(X_{\infty})),Z/l^{\nu}) \mapsto 0\]

  \vspace{30pt}

 Under our assumptions, the left hand side would be $0$ since multiplying that infinite product of quasicyclic groups by a power of $l$ would be equal to the initial infinite product of quasicyclic $l$ groups because of the divisible structure that the $l$ quasicyclic groups have.Therefore since the middle term in the above exact sequence which is $\pi_{m}(LK(X_{\infty}))\wedge M(Z/l^{\nu})$ is finite by page 387 in [7]and page 388 (13) in [7]  it would be isomorphic to the $l^{\nu}$ torsion elements of  $\pi_{m}(LK(X_{\infty}))$ Then this torsion group would have to be finite and therefore it cannot have as a structure an infinite product of quasicyclic groups since that group has an $l^{\nu}$ torsion subgroup which is infinite because of the structure the group has. We can conclude then that case a) is not possible.

 \vspace{15pt}

   Now we turn to case b)in which  $\pi_{m}((LK(X_{\infty})))$ has as a subgroup included in $A_{r}$ an infinite product of finite cyclic groups $Z/l^{\nu}Z$ with $\nu$ growing to infinity or not growing to infinity. If $\nu$ does not grow to infinity the proof is exactly as the proof of case a). If $\nu$ grows to infinity We have
  $L=\Pi_{\nu=1}^{\infty} Z/l^{\nu}Z$ as a subgroup of $A_{r}$, so that it is sufficient to consider $A_{r}=L\times H$ and then $\pi_{m}(LK(X_{\infty}))=L\times H\times A_{d}$. Replacing in the exact splitting sequence (1) we get

\vspace{30pt}
\[0 \mapsto L/l^{\nu}L\times H/l^{\nu}H\times A_{d}/l^{\nu}A_{d} \mapsto  \pi_{m}(LK(X_{\infty})/l^{\nu}) \mapsto \tag{6}\]

\vspace{3pt}

\[\mapsto Tor^{1} (\pi_{m-1}(LK(X_{\infty})),Z/l^{\nu}) \mapsto 0\]

\vspace{30pt}

 But observe that $L/l^{\nu}L$ is isomorphic to $\Pi_{j=1}^{j=\nu} Z/l^{j}Z\times \Pi_{j=\nu+1}^{\infty} Z/l^{\nu}Z$

\vspace{5pt}

 which is an infinite abelian group and that cannot be possible since the splitting exact sequence (6) has a finite abelian group as middle term as already pointed out above and in [7] pp 387 and pp 388 (13).

\vspace{30pt}

Therefore neither case a) or b) can hold and $Hom(Z/l^{\nu}Z,\pi_{m}(LK(X_{\infty})))$ must be at most countable as wanted.

\vspace{30pt}

Remark:2.4 : this theorem 2.3 shows what we conjectured above in order to research the Tate module M, that homotopy groups of the $l$-adic completion of the $LK(X_{\infty})$ spectrum are isomorphic to the $l$-adic completion of the homotopy groups of the K spectrum, which at first sight is not trivial.

 \vspace{25pt}

Corollary 2.5: \emph{$\pi_{-1}(LK(X_{\infty}))$ is reduced}

\vspace{100pt}

\printbibliography % Prints the bibliography

References

\vspace{25pt}

[1] Ioannis Emmannouil
\emph{Mittag-Leffler Condition and the vanishing of the derived inverse limit}
Topology Vol 35 No 1 pp 267-271  1996

\vspace{8pt}

[2] B Gray
\emph{Spaces of the same n-type for all n}
Topology 5 (1966)  pp 241-243

\vspace{8pt}

[3] R Jardine
\emph{Generalized Etale Cohomology Theories}
Progress in Mathematics Vol 146 Birkhauser-Verlag

\vspace{8pt}

[4]S.A.Mitchell
\emph{Hypercohomology Spectra andThomason Descent Theorem}
Algebraic K Theory 16  221-277

\vspace{8pt}

[5]Roitman
\emph{An Introduction to Homological Algebra}Academic Press

\vspace{8pt}

[6]R.W.Thomason
\emph{Algebraic K-Theory and Etale Cohomology}
Ann.Ec.Norm.Sup 437-552

\vspace{8pt}

[7]R.W.Thomason.
\emph{A Finiteness Condition Equivalent to the Tate Conjecture over Fq}
Contemporary Mathematics \textbf{83} 385-392.

\end{document}